\documentclass[11pt]{article}

\usepackage[applemac]{inputenc} 
\usepackage{amsthm,amssymb,amsbsy,amsmath,amsfonts,amssymb,amscd,mathrsfs}

\long\def\symbolfootnote[#1]#2{\begingroup%
\def\thefootnote{\fnsymbol{footnote}}\footnote[#1]{#2}\endgroup} 

\usepackage{graphicx}
\DeclareGraphicsExtensions{.jpg}
\usepackage{color}
\usepackage{fullpage}
\usepackage{ifsym}
\usepackage[hypertex]{hyperref}
\newtheorem{theorem}{Theorem}
\newtheorem{corollary}[theorem]{Corollary}
\newtheorem{lemma}[theorem]{Lemma}
\newtheorem{proposition}[theorem]{Proposition}
\theoremstyle{definition} 
\newtheorem{definition}[theorem]{Definition}
\newtheorem{example}[theorem]{Example}
\theoremstyle{remark}
\newtheorem{remark}[theorem]{Remark}
\newcommand{\bt}{\begin{theorem}}
\newcommand{\et}{\end{theorem}}
\newcommand{\bl}{\begin{lemma}}
\newcommand{\el}{\end{lemma}}
\newcommand{\bp}{\begin{proposition}}
\newcommand{\ep}{\end{proposition}}
\newcommand{\bc}{\begin{corollary}}
\newcommand{\ec}{\end{corollary}}
\newcommand{\bdeff}{\begin{definition}}
\newcommand{\edeff}{\end{definition}}
\newcommand{\brem}{\begin{remark}}
\newcommand{\erem}{\end{remark}}
\newcommand{\bex}{\begin{example}}
\newcommand{\eex}{\end{example}}

\renewcommand{\r}[1]{(\ref{#1})}
\newcommand{\con}{{\mathcal C}}
\newcommand{\mc}[1]{\mathcal {#1}}

\newcommand{\bi}{\begin{itemize}}
\newcommand{\iii}{\item}
\newcommand{\ei}{\end{itemize}}
\newcommand{\bd}{\begin{description}}
\newcommand{\ed}{\end{description}}
\newcommand{\all}{\forall \,}
\newcommand{\bqn}{\begin{equation}}   
\newcommand{\eqn}{\end{equation}}       
\newcommand{\bqna}{\begin{eqnarray}}
\newcommand{\eqna}{\end{eqnarray}}
\newcommand{\eqnan}{\nonumber\end{eqnarray}}

\newcommand{\nn}{\nonumber}
\newcommand{\ba}[1]{\begin{array}{#1}}
\newcommand{\ea}{\end{array}}

\newcommand{\R}{\mathbb{R}}

\newcommand{\virg}[1]{``#1''}

\newcommand{\lam}{\lambda}
\newcommand{\g}{\gamma}

\newcommand{\eps}{\varepsilon}

\newcommand{\tx}[1]{\mathrm{#1}}
\newcommand{\til}[1]{\widetilde{#1}}
\newcommand{\la}{\left\langle}
\newcommand{\ra}{\right\rangle}

\newcommand{\distr}{\Delta}
\newcommand{\metr}{{\bf g}}

\newcommand{\Pg}[1]{\left\{ #1 \right\}}

\newcommand{\Exp}{\mc{E}}

\newcommand{\Q}{Q}
\newcommand{\Qh}{\widehat{Q}}
\newcommand{\Bad}{\Sigma}
\newcommand{\Badinf}{\Sigma_{\infty}}
\newcommand{\Badfin}{\Sigma_{0}}
\newcommand{\wh}{\widehat}
\newcommand{\A}{A}
\newcommand{\z}{\alpha}
\newcommand{\tr}{\text{trace}}
\newcommand{\tcut}{t_{cut}}
\newcommand{\tconj}{t_{conj}}
\newcommand{\m}{\frac{1}{2}}
\newcommand{\h}{\mathfrak{h}}
\newcommand{\lv}{\mathfrak{L}}

\begin{document}
\begin{center} \noindent
{\LARGE{\sl{\bf On 2-step, corank 2 nilpotent sub-Riemannian metrics}}}
\vskip 0.6 cm
Davide Barilari\\ 
{\footnotesize SISSA, Via Bonomea 265, Trieste, Italy - {\tt barilari@sissa.it}}\\
\vskip 0.3cm
Ugo Boscain\symbolfootnote[0]{This research has been supported by the European Research Council, ERC StG 2009 \virg{GeCoMethods}, contract number 239748, by the ANR Project GCM, program \virg{Blanche}, project number NT09-504490 and by the DIGITEO project CONGEO.}\\
{\footnotesize CNRS, CMAP Ecole Polytechnique and equipe INRIA GECO, Paris, France - {\tt boscain@cmap.polytechnique.fr}}
\vskip 0.3cm
Jean-Paul Gauthier\\ 
{\footnotesize Laboratoire LSIS, Université de Toulon, and equipe INRIA GECO, France - {\tt gauthier@univ-tln.fr}}\\

\vskip 0.3cm
\end{center}
\vskip 0.1 cm

\begin{center}
\today
\end{center}
\vskip 0.2 cm
\begin{abstract} 

In this paper we study the nilpotent 2-step, corank 2 sub-Riemannian metrics that are nilpotent approximations of general sub-Riemannian metrics. We exhibit optimal syntheses for these problems. It turns out that in general the cut time is not equal to the first conjugate time but has a simple explicit expression. As a byproduct of this study we get some smoothness properties of the spherical Hausdorff measure in the case of a generic 6 dimensional,  2-step corank 2 sub-Riemannian metric.

\end{abstract}
 
\section{Introduction}
\subsection{Sub-Riemannian manifolds}\label{s:pmp}
In this paper, by a sub-Riemannian manifold we mean a triple $(M,\distr,\metr)$, where $M$ is a connected orientable smooth manifold of dimension $n$, $\distr$ is a smooth vector distribution of constant rank $m<n$, satisfying the H\"ormander condition and $\metr$ is an Euclidean structure over $\distr$. 

A sub-Riemannian manifold has the natural structure of a metric space, where the distance is the so called Carnot-Caratheodory distance
\begin{eqnarray}
d(q_0,q_1)=
\inf\{\int_0^T\sqrt{\metr_{\g(t)}(\dot\g(t),\dot\g(t))}~dt~|~ \g:[0,T]\to M \mbox{ is a Lipschitz curve},\cr
\g(0)=q_0,\g(T)=q_1, ~~\dot \g(t)\in\distr_{\g(t)}\mbox{ a.e. in $[0,T]$} \}.\nonumber
\end{eqnarray}
As a consequence of the H\"ormander condition $d$ is actually a distance inducing the topology of the manifold. This is the  Rashevsky-Chow Theorem,  see for instance \cite{agrachevbook} or \cite{gromov}.

Along this paper we assume that the structure is 2-step bracket generating i.e. 
$$
T_qM=\Delta_q+[\Delta,\Delta]_q,~~\mbox{ for every }q\in M,
$$
and we quote a 2-step sub-Riemannian metric by its rank and its dimension, i.e. with the pair $(m,n)$. The quantity $k=n-m$ is called the {\it corank} of the structure.

It is know from Mitchell \cite{mitchell} that the Hausdorff dimension of $M$, as a metric space, is  $\h=2n-m>n$. In this paper we focus on the case $(m,m+2)$.

A sub-Riemannian manifold is left-invariant if $M=G$, a Lie group, and both $\distr$ and $\metr$ are left-invariant over $G$.

Locally, the pair $(\distr,{\mathbf g})$ can be specified by the data of a set of $m$ smooth vector fields spanning $\distr$, being an orthonormal frame for ${\mathbf g}$, i.e.  
\bqn
\label{trivializable}
\distr_{q}=\text{span}\{X_1(q),\dots,X_m(q)\}, \qquad \qquad \metr_q(X_i(q),X_j(q))=\delta_{ij}.
\eqn
In this case, the set $\Pg{X_1,\ldots,X_m}$ is called a \emph{local orthonormal frame} for the sub-Riemannian metric. 

The sub-Riemannian metric can also be expressed locally in ``control form'' as follows. We consider the control system,
\bqn
\dot q=\sum_{i=1}^m u_i X_i(q)\,,~~~u_i\in\R\,,
\eqn
and the problem of finding the shortest curve minimizing that joins two fixed points $q_0,~q_1\in M$ is naturally formulated as the optimal control problem,
\bqn \label{eq:lnonce}
~~~\int_0^T  \sqrt{
\sum_{i=1}^m u_i^2(t)}~dt\to\min, \ \ \qquad q(0)=q_0,~~~q(T)=q_1.
\eqn
A {\it geodesic} for the sub-Riemannian metric is a curve, parametrized by constant velocity, such that every short enough piece of it is a local minimizer of the length.  For the sub-Riemannian metrics given in control form all the geodesics can be computed with   Pontryagin's maximum principle \cite{pontryaginbook}. In the 2-step bracket generating case, it is known that there is no strict abnormal minimizer, and all geodesics are projections on $M$ of the trajectories of the Hamiltonian system associated with the following Hamiltonian over $T^\ast M$
\bqn \label{eq:ham00}
H(\lam,q)=\frac12 \sum_{i=1}^m\langle \lam,X_i(q)\rangle^2, \qquad \lam \in T^{*}_{q}M.
\eqn
and corresponding to the level set $\{H=c\}$, for $c>0$.

\subsection{Nilpotent approximation}

Consider a sub-Riemannian manifold $(M,\distr,\metr)$ and fix a point $q\in M$. The Lie bracket induces a skew symmetric tensor bilinear mapping
\bqn \label{eq:bil}
[\cdot,\cdot]_{q}: \Delta_q\times\Delta_q\to T_qM/\Delta_q.
\eqn
Then, for every $Z^\ast\in(T_qM/\Delta_q)^\ast$, we have
$$
Z^\ast([X,Y]+\Delta_q)=\langle A_{Z^\ast}(X),Y\rangle_\metr,
$$
for some $\metr$-skew symmetric endomorphism $A_{Z^\ast}$ of $\Delta_q$.

\brem[Notation]\label{r-LFQ}
We denote by $\mc{L}_q$ the $k$-dimensional space of skew symmetric endomorphisms of $\Delta_q$ obtained by taking the union of all the $A_{Z^\ast}$ at $q$. This notation is used in the Appendix.
\erem

The space $L_{q}=\Delta_q\oplus T_q M/\Delta_q$ is endowed with  the structure of a 2-step nilpotent Lie-algebra by setting 
$$
[(V_{1},W_{1}),(V_{2},W_{2})]=(0,[V_{1},V_{2}]+\Delta_q).
$$ 
The associated simply connected nilpotent Lie group is denoted by $G_q$ and the exponential mapping Exp$:L_q\to G_q$ is one to one and onto. By translation, the metric $\metr_q$ over $\Delta_q$ allows to define a left-invariant sub-Riemannian metric over $G_q$.
\bdeff \label{d:richiamo} 
The sub-Riemannian metric on $G_{q}$ defined above is called the \emph{nilpotent approximation} of $(M,\Delta,\metr)$ at $q$. 
\edeff
Any $k$ dimensional vector sub-space ${\cal V}_q$ of $T_qM$, transversal to $\Delta_q$ allows to identify $L_q$ and $G_q$ to $T_q M\simeq\Delta_q\oplus T_qM/\Delta_q$.

Fix $q_0\in M$. We  can chose coordinates $x$ in $\Delta_{q_0}$ such that the metric $\metr_{q_0}$ is the standard Euclidean metric, and for any linear coordinate system $y$ in ${\cal V}_{q_0}$, there are skew symmetric matrices $L_1,\ldots,L_k\in$ $so(m)$ such that the mapping \eqref{eq:bil} writes
$$
[X,Y]+\Delta_{q_0}= \left(\ba{c} X' L_1 Y\\ \vdots \\ X' L_{k} Y\ea\right).
$$
where $X^{\prime}$ denotes the transpose of the vector $X$.
Then the nilpotent approximation written in control form is

\bqn \label{eq:system0}
\begin{cases}
\dot{x}_{i}=u_{i}, \qquad i=1,\ldots,m,\\
\dot{y}_1=\m x'L_1 u,\\
\quad\vdots\\
\dot{y}_{k}=\m x'L_{k} u.
\end{cases}
\eqn

The construction of the nilpotent approximation given in Definition \ref{d:richiamo} makes sense for any sub-Riemannian metric, but it coincides with the standard one (see \cite{nostrolibro,bellaiche}) in the 2-step bracket generating case only. 

\bp
The distribution is 2-step bracket generating if and only if the endomorphisms of $\Delta_q$, $L_i$, $i=1,\ldots,k$ (respectively the matrices $L_i$ when coordinates $y$ in ${\cal V}_q$ are chosen) are independent.
\ep
In the 2-step bracket generating case these linear coordinates $y$ in $T_qM/\Delta_q$ may be chosen in such a way that the endomorphisms $L_i$, $i=1,\ldots,k$ are orthonormal with respect to the Hilbert-Schmidt norm $\langle L_i,L_j\rangle=\frac{1}{m}\text{trace}(L_i' L_j)$. This choice defines a canonical Euclidean structure in $T_qM/\Delta_q$ and a corresponding volume in $T_qM/\Delta_q$. Then using the Euclidean structure over $\Delta_q$ we get a canonical Euclidean structure over $\Delta_q\oplus T_qM/\Delta_q$. The choice of the vector subspace ${\cal V}_q$ induces an Euclidean structure on $T_qM$ which depends on the choice of ${\cal V}_q$, but the associated volume on  $T_qM$ is independent on this choice. 
\bdeff
This volume form on $M$ is called the Popp measure.
\edeff
\noindent
The Popp measure is a smooth volume form.

\subsection{Statement of the results}

\subsubsection{History}

The main purpose of  this paper is to  build the optimal synthesis for $(m,m+2)$ nilpotent sub-Riemannian metrics, i.e. the set of all trajectories starting from the identity of the group and realizing the minimum of the distance, with a precise description of their cut time.

Optimal syntheses are in general very difficult to obtain.  Usually the steps are the following:
\bi
\iii[-] Apply first order necessary conditions for optimality (which in the case of sub-Riemannian manifolds are given by the Pontryagin Maximum Principle) to reduce the set of candidate optimal trajectories. This first step can be already very difficult since one should find solutions of a Hamiltonian system, which  is not integrable in general.
\iii[-] Use higher order necessary conditions to reduce  further the set of optimal trajectories. This step  usually leads to the  computation of  the conjugate locus, i.e. the set of points up to which geodesics are locally optimal.
\iii[-] Prove that no strict abnormal extremal is optimal (for instance by using conditions such as the so called  Goh condition \cite{nostrolibro,agrachevbook}).  If one fails to go beyond this step, then one can hardly get an optimal synthesis, since no general technique exists to treat abnormal minimizers.  
\iii[-]  Among all solutions of the first order necessary conditions, find the optimal ones. One has to  prove that, for each point of a candidate optimal trajectory, there is no other trajectory among the selected ones, reaching that point. The first point after which a first order trajectory loses global optimality  is called a {\it cut point}. The union of all cut points is the {\it cut locus}.
\ei
As a consequence of these difficulties, optimal syntheses in sub-Riemannian geometry have been obtained in few cases.

The most studied cases are those of left invariant sub-Riemannian metrics (see for instance \cite{miosr3d} for a classification in the 3D case). The first optimal synthesis was obtained for the Heisenberg group in \cite{gaveau,gershkovich}.
Then complete optimal syntheses were obtained for the 3D simple Lie groups $SU(2)$, $SO(3)$, $SL(2)$, with the metric induced by the Killing form in \cite{francesco-SIAM,boscainrossiRSC}. An impressive work has been done by Yuri Sachkov who obtained the optimal synthesis for the group of motions of the plane $SE(2)$ (see \cite{yuri1,yuri2}).

In dimension larger than 3, only nilpotent groups have been attacked. The complete  optimal synthesis was obtained in \cite{corank1} in the contact nilpotent case. Some results were obtained by Y. Sachkov for the Engel  and Cartan groups \cite{yuri234,yuri235}.

When a Lie group structure is not available there are also some results: the optimal synthesis was obtained for a neighborhood of the starting point in the 3D contact case in \cite{agrachev-contact,charlotpazzesco,gauthier-contact} and in the 4D quasi-contact case in \cite{gregoire}. The optimal synthesis was obtained in the important Martinet nilpotent case,  where  abnormal minimizers can be optimal (see \cite{bonnard-kupka}). They also solved the problem for certain perturbations of this case where strictly abnormal minimizers occur (see \cite{bonnard-trelat}).

To our knowledge, no other case has been solved.

It is interesting to notice that when the sub-Riemannian metric is invariant by certain continuous transformation (e.g. rotations) then most of the cut points are  automatically conjugate. This happens for instance on the Heisenberg group, on $SU(2)$ and in the contact nilpotent case.

\brem[Notation] 
In the case of our nilpotent approximations, covectors in $T^\ast_qM$ can be identified with vectors in $T_qM$ via the Euclidean structure of $T_qM$ given by the choice of ${\cal V}_q$. In our coordinates $(x,y)$, these covectors/vectors are typically denoted by $(u_0,r)$.   
\erem

For nilpotent $(m,m+1)$ sub-Riemannian metrics that are nilpotent approximations of general sub-Riemannian metrics, the control  systems can be written as
 
\bqn \label{s:1} \tag{C1}
\begin{cases}
\dot{x}_{i}=u_{i}, \qquad i=1,\ldots,m,\\
\dot{y}=\m x^{\prime}L u,\qquad L \ \text{skew symmetric,}
\end{cases}
\eqn
Denote by $\sigma(L)$ the set of all moduli of eigenvalues of the matrix $L$.
The following fact is proved in \cite{corank1}:
\bt (\cite{corank1})
Arclength geodesics of system \eqref{s:1}, starting from the origin, are parametrized by an initial covector $\lam_{0}=(u_{0},r)\in S^{m-1}\times  \R$, and they 
are optimal until time
$$t_{cut}(\lam_{0})=\frac{2\pi}{|r|  \max \sigma(L)},$$ 
 with the understanding $\tcut(\lam_{0})=+\infty$ if $r=0$.
Moreover $t_{cut}(\lam_{0})=t_{conj}(\lam_{0})$.
\et
The proof of this result is based on the fact that geodesics can be expressed in terms of usual trigonometric functions and, thanks to a certain monotonicity property, the cut locus can be explicitly computed and is exactly equal to the conjugate locus. 

\subsubsection{Optimal synthesis for the nilpotent $(m,m+2)$ case}

The main result of this paper is the optimal synthesis in the case of a nilpotent approximation in the $(m,m+2)$ case. In this case  the control system  
can be written in coordinates 
$q=(x_{1},\ldots,x_{k},y_{1},y_{2})$ as
\bqn \label{s:2} \tag{C2}
\begin{cases}
\dot{x}_{i}=u_{i}, \qquad i=1,\ldots,m,\\
\dot{y}_{1}=\m x^{\prime}L_{1}u, \\
\dot{y}_{2}=\m x^{\prime}L_{2} u,\qquad L_{1},L_{2} \ \text{skew symmetric}.
\end{cases}
\eqn
   Set $r_1=|r|\cos\theta,$ $r_2=|r|\sin\theta,$ and $L_{\theta}=\cos (\theta) L_{1}+\sin (\theta) L_{2}$.
 \bt \label{t:ttt}
   \label{t-synthesis} 
Arclength geodesics of system \eqref{s:2}, starting from the origin, are parametrized by an initial covector $\lam_{0}=(u_{0},r)\in S^{m-1}\times  \R^{2}$, and they 
are optimal until time
$$\tcut(\lam_{0})=\frac{2\pi}{ \max \sigma(r_{1}L_{1}+r_{2}L_{2})}=\frac{2\pi}{|r|  \max \sigma(L_{\theta})},$$ with the understanding $\tcut(\lam_{0})=+\infty$ if $r=0$.  Moreover, in general, $t_{cut}(\lam_{0})\neq t_{conj}(\lam_{0})$.
\et  

The last statement in the Theorem \ref{t:ttt} says that in the corank 2 case the cut and the conjugate time coincide only in some particular cases, which we explicitly describe in the $(4,6)$ case (see Theorem \ref{t-quaternions}).

The reason why the corank 2 case is more difficult than the corank 1 case is precisely the fact that the cut locus is not equal to the conjugate locus. (The latter we are not able to compute explicitly.)

Explicit expression of geodesics for this optimal synthesis are given in Section \ref{s:cutconj}.

\subsubsection{The nilpotent $(4,6)$ case}
In the nilpotent $(4,6)$ case our first result is the following:
\bt
\label{t-quaternions}
The following properties are equivalent:
\bi
\iii[(P1)] The first  conjugate locus is equal to the cut locus.
\iii[(P2)] The linear coordinates $y$ in $T_qM/\Delta_q$ can be chosen in such a way that the pair  $(L_1,L_2)$ of $4\times4$ skew symmetric matrices belongs to the set $({\Q\cup \Qh})^2$. 
\ei
 \et
Here $\Q$ (resp. $\Qh$) denotes the set of pure quaternions (resp. pure skew quaternions), see Appendix \ref{ap-quaternions}.

Our second result is a continuation of the paper \cite{corank1} for corank 1, where the following result is proved for general sub-Riemannian metrics.

\bt[\cite{corank1}]  In the $(m,m+1)$ case the Radon-Nykodym derivative of the spherical Hausdorff measure with respect to the Popp measure is a ${\cal C}^3$ function, but is not ${\cal C}^5$ in general. 
\et
Here we show the following result
\bt
\label{t-haus}
For  a generic $(4,6)$ sub-Riemannian metric\footnote{which means for an open and dense subset of all $(4,6)$ sub-Riemannian metrics, endowed with the Whitney topology.}, the Radon-Nykodym derivative of the spherical Hausdorff measure with respect to the Popp measure is ${\cal C}^1$.
\et
In the previous paper \cite{corank1} it is shown that the Radon-Nikodym derivative of the spherical Hausdorff measure with respect to the Popp measure is inversely proportional (as a function of $q$) to the volume of the unit sub-Riemannian ball of the nilpotent approximation at $q$. Then Theorem \ref{t-haus} is a byproduct of the optimal synthesis given here.

Note that in the corank 1 case, the higher differentiability of the Radon-Nikodym derivative is due to the fact that the conjugate locus is equal to the cut locus, which is not the case here. 

Due to the complexity of the computations even in this low dimensional case, it is not easy to determine the real degree of differentiability of Hausdorff measure. This is still an interesting open question.

\subsection{Organization of the paper}

Section \ref{s:cutconj} is devoted to the construction of the optimal synthesis for $(m,m+2)$ nilpotent sub-Riemannian metrics and, as a consequence, to the proof of Theorem \ref{t-synthesis}. In Sections \ref{s-exp} and \ref{s:exp2} we compute the exponential map. In Section \ref{s:cut} we prove that geodesics are optimal up to $\tcut$. Finally in Section \ref{s:conj} we show that the cut time does not coincide, in general, with the first conjugate time. 
In Section \ref{s-46} we give the proofs of Theorems \ref{t-quaternions} and \ref{t-haus}.

In the Appendix we recall basic facts about quaternions, we prove a technical Lemma, and applying an Abraham's transversality theorem, we prove that, generically, for the $(4,6)$ case, a certain \virg{bad set} is made of isolated points, which permits to conclude about the differentiability of the Radon-Nikodym derivative (Theorem \ref{t-haus}).

\section{Exponential map and synthesis}\label{s:cutconj}

\subsection{Hamiltonian equations in the $(m,n)$ case}\label{s-exp}
The purpose of this section is to compute the exponential map, i.e. the set of all geodesics, parametrized by length, starting form the origin of the control system \eqref{eq:system0}, i.e. the system
\bqn \label{eq:controlsyst}
\begin{cases}
\dot{x}_{i}=u_{i}, \qquad\qquad i=1,\ldots,m,\\
\dot{y}_{h}=\m x^{\prime}L_{h}u, \qquad h=1,\ldots,k.
\end{cases}
\eqn
Let $L_{h}=(b_{ij}^{h}),$ for $h=1,\ldots,k$. Then the control system can be written in the form
$\dot{q}=\sum_{i=1}^{m}u_{i}X_{i}(q)$ where $q=(x,y)$ and 
\begin{gather}
X_i=\partial_{x_i}+\m \sum_{j,h} b_{ij}^{h} x_j \partial_{y_{h}}, \qquad i=1,\ldots,m. \nn
\end{gather}
Setting $Y_{h}=\partial_{y_{h}},$ for  $h=1,\ldots,k,$ the commutation relations are
\begin{gather}
[X_i,X_j]=\sum_{h=1}^{k}b_{ij}^h Y_h,  \qquad i,j=1,\ldots,m, \label{eq:xixj}\\
[X_i,Y_j]=[Y_j,Y_h]=0, \qquad i=1,\ldots,m, \qquad j,h=1,\ldots,k. \label{eq:xiyj}
\end{gather}
Define the functions on $T^{*}M$, that are linear on fibers,
\begin{gather}
u_{i}(\lam,q)=\la \lam, X_{i}(q)\ra, \qquad i=1,\ldots,m,\\
r_{h}(\lam,q)=\la \lam, Y_{h}(q)\ra, \qquad h=1,\ldots,k.
\end{gather}
These functions can be treated as coordinates on the fibers of $T^{*}M$ to solve the Hamiltonian system given by the Pontryagin Maximum Principle, see Section \ref{s:pmp}. This Hamiltonian system is associated with the Hamiltonian 
\bqn\label{eq:Ham1}
H(\lam,q)=\frac12 \sum_{i=1}^m\langle \lam,X_i(q)\rangle^2= \m \sum_{i=1}^{m}u_{i}^{2}(\lam,q), \qquad \lam \in T^{*}_{q}M.
\eqn
\brem \label{r:time} The geodesics parametrized by length correspond to the level set $\{H=1/2\}$. Notice that, for systems of type $\dot{q}=\sum_{i=1}^{m}u_{i}X_{i}(q)$, with fixed initial and final points, the problem of finding length-parametrized curves minimizing the length, is equivalent to the problem of minimizing time with the constraint $\{\|u\|\leq 1\}$.\erem
For a function $a\in C^{\infty}(T^{*}M)$ we have that, along the sub-Riemannian flow
\bqn \label{eq:derh}
\dot{a}=\{a,H\}=\sum_{i=1}^{m} \{a,u_{i}\}u_{i},
\eqn
where $\{a,b\}$ denotes the Poisson bracket of two functions in $T^{*}M$. The following Lemma gives a way of computing the covector $\lam(t)$, solution of the Hamiltonian system associated with \eqref{eq:Ham1} in the coordinates $(u,r)$. 
\bl
If $u(t)$ and $r(t)$ are solution of \eqref{eq:derh} corresponding to level set $\{H=1/2\}$, then they satisfy 
$$\begin{cases}
\dot{u}(t)=(r_{1}L_{1}+\ldots+r_{k}L_{k}) u(t), \qquad u(0)=u_{0}, \qquad \|u_{0}\|=1,\\
\dot{r}(t)=0.
\end{cases}
$$
\el
\begin{proof}Remind that, if $a_{i}(\lam,q)=\la \lam,Z_{i}(q)\ra$, for some vector fields $Z_{i}$, $i=1,2$, then
$$\{a_{1},a_{2}\}=\la\lam,[Z_{1},Z_{2}]\ra.$$ 
Applying \eqref{eq:derh} for $a=r_{h}$ and using \eqref{eq:xiyj} we get
$$\dot{r}_{h}=\sum_{h=1}^{k} \{r_{h},u_{i}\}u_{i}=0 \quad \Rightarrow \quad  r_{h}=\tx{const}.$$
Similarly, using \eqref{eq:xixj}, one find 
\begin{gather*}
\dot{u}_{i}=\sum_{i=1}^{n} \{u_{i},u_{j}\}u_{j}=\sum_{i=1}^{n} b_{ij}^{k}r_{k} u_{j}.
\end{gather*}
\end{proof}
\brem
In the following geodesics are parametrized by the initial covector $\lam(0)=(p(0),r(0))=(u_{0},r)$, since $r=const$ and $u_{i}=\la \lam, X_{i}\ra$ and at the starting point we have $X_{i}(0)=\partial_{x_{i}}$.
\erem
\subsection{Exponential map in the corank 2 case} \label{s:exp2}
From now on we focus on the case $(m,m+2)$, i.e. when the corank $k$ is equal to 2.
We can write the equation of geodesics starting from the origin as follows
\bqn \label{eq:qqq0}
\begin{cases}
\displaystyle{x(t)=\int_{0}^{t}e^{s(r_{1}L_{1}+r_{2}L_{2})}u_{0}ds} , \quad x(0)=0,\\[0,1cm]
\displaystyle{y_{1}(t)=\m \int_{0}^{t}x(s)^{\prime} L_{1}\, u(s) ds}, \quad y_{1}(0)=0,\\[0,1cm]
\displaystyle{y_{2}(t)=\m \int_{0}^{t}x(s)^{\prime} L_{2} \,u(s) ds}, \quad y_{2}(0)=0.
\end{cases}
\eqn

\brem[Notation] In the following we denote by $E^{u_{0},r_{1},r_{2}}_{L_{1},L_{2}}(t)$ the geodesic, parametrized by the length, and starting from the origin, defined by equations \eqref{eq:qqq0}, associated with $L_{1},L_{2}$.
\erem
\bdeff
The matrices $L_{1},L_{2}$ being fixed, the \emph{exponential map} is the map $\Exp: \R^{+} \times \Lambda  \to \R^{n}$ defined by 
$$\Exp(t,u_{0},r_{1},r_{2})= E_{L_{1},L_{2}}^{u_{0},r_{1},r_{2}}(t), \qquad \Lambda=\{(u_{0},r_{1},r_{2}), u_{0}\in S^{m-1}, r_{i}\in \R\}$$  
\edeff

\brem
The optimal control problem 
\bqn \label{eq:controlsyst00}
\begin{cases}
\dot{x}_{i}=u_{i}, \qquad\qquad i=1,\ldots,m,\\
\dot{y}_{1}=\m x^{\prime}L_{1}u,\\
\dot{y}_{2}=\m x^{\prime}L_{2}u,
\end{cases}
\eqn
 is invariant with respect to the following change of coordinates
\bi
\iii[(a)] orthogonal changes of coordinates in the $x$ space,
\iii[(b)] linear changes of coordinates in the $y$ space.
\ei
Indeed, let $M$ be a nonsingular orthogonal matrix $(M^{-1}=M^{\prime})$ and define the new coordinates
$
\til{x}=Mx$. Then 
$$\dot{\til{x}}=M\dot{x}=Mu=:\til{u},$$
and
$$\dot{y}_{i}=x^{\prime}L_{i} u=(Mx)^{\prime} M L_{i}M^{\prime} (Mu)=\til{x}^{\prime} M L_{i}M^{\prime} \til{u}.$$
Hence, in the new coordinates, 
$L_{i}$ is changed for $\til{L}_{i}:=M L_{i}M^{\prime}.$ 

Also, it is easy to see that the change of coordinates 
\bqn\label{eq:lcc}
\til{y}_{1}=\alpha_{1}y_{1}+\alpha_{2}y_{2}, \qquad \til{y}_{2}= \beta_{1}y_{1}+\beta_{2}y_{2},
\eqn
corresponds to the change 
$$\til{L}_{1}=\alpha_{1}L_{1}+\alpha_{2}L_{2}, \qquad \til{L}_{2}= \beta_{1}L_{1}+\beta_{2}L_{2}.$$
In other words we can change $L_{1}$ and $L_{2}$ up to congruence and linear combinations.

\erem
Using these arguments one immediately gets
\bl \label{l:reducer0}
Let $(r_{1},r_{2})=:(r\cos \theta,r \sin \theta)$ and $L_{\theta}:=\cos \theta L_{1}+\sin \theta L_{2}$, $\til{L}_{\theta}:=-\sin \theta L_{1}+\cos \theta L_{2}$. Consider the rotation matrix
$R_{\theta}=\left(
\begin{smallmatrix}
\cos \theta&\sin \theta\\-\sin \theta& \cos \theta
\end{smallmatrix}
\right)
$
and the orthogonal matrix $M$ such that $ML_{\theta}M'$ is block diagonal.

Denote $\Omega=\left(
\begin{smallmatrix}
M&0\\0&R_{\theta}
\end{smallmatrix}
\right)$ and $\til{u}_{0}=Mu_{0}$.
We have the equality
\bqn \label{eq:EE}
 \Omega\, E^{u_{0},r_{1},r_{2}}_{L_{1},L_{2}}(t)=
E^{\til{u}_{0},r,0}_{L_{\theta},\til{L}_{\theta}}(t).
\eqn
\el
\medskip
Thanks to Lemma \ref{l:reducer0}, one can always restrict to geodesics of the type

\bqn \label{eq:qqq1}
\begin{cases}
\displaystyle{x(t)=\int_{0}^{t}e^{srL_{\theta}}u_{0}ds}, \\[0,1cm]
\displaystyle{y_{1}(t)=\m \int_{0}^{t}x(s)^{\prime} L_{\theta}\, u(s) ds},\\ 
\displaystyle{y_{2}(t)=\m \int_{0}^{t}x(s)^{\prime} \til{L}_{\theta} \,u(s) ds},
\end{cases}
\eqn
where
 $L_{\theta}$ is in the block-diagonal form
$$L_{\theta}=
\begin{pmatrix}
0&a_{1}&&&\\
-a_{1}&0&&&\\
&&\ddots&&\\
&&&0&a_{\ell}\\
&&&-a_{\ell}&0
\end{pmatrix}, \qquad \text{or} \qquad \begin{pmatrix}
0&a_{1}&&&&\\
-a_{1}&0&&&&\\
&&\ddots&&&\\
&&&0&a_{\ell}&\\
&&&-a_{\ell}&0&\\
&&&&&&0\\
\end{pmatrix},$$
depending on the fact that $m$ is even ($m=2\ell$) or odd ($m=2\ell+1$), and where the
geodesic is associated with the covector $(r_{1},r_{2})=(r,0)$.

\brem
When we deal with a fixed sub-Riemannian metric we can assume also that the coordinates in the $x$ space are chosen in such a way that $a_{1} \geq a_{i}$, for every $i$. In this case
$$\frac{2\pi}{r  \max (\sigma(L_{\theta}))}=\frac{2\pi}{a_{1}r}.
$$
\erem



\subsection{Computation of the cut time} \label{s:cut}
In this section we prove Theorem \ref{t-synthesis}, i.e. we compute the last time at which a geodesic parametrized by length  is optimal.

We first consider the case $r=0$. In this case equations \eqref{eq:qqq0} can be easily integrated and gives the straight lines
$$
\begin{cases}
x(t)=u_{0}t,\\
y_{i}(t)=0.
\end{cases}
$$
This trajectory is optimal for any time (i.e. $\tcut=+\infty$) since the sub-Riemannian length of a geodesic coincides with the Euclidean length of its projection on the horizontal subspace $(x_{1},\ldots,x_{m})$, as follows from formula \eqref{eq:lnonce}.
\\

In what follows we use the notation $A=rL_{\theta}$, $\til{A}=r\til{L}_{\theta}$ and we focus on the case when $A$ is even dimensional (i.e. $m=2\ell$) and invertible. The case $A$ non invertible (in particular $A$ odd dimensional) needs an obvious modification of the proof. 

With this notation the system \eqref{eq:qqq1} is rewritten as
\bqn \label{eq:qqq3}
\begin{cases}\displaystyle{
x(t)=A^{-1}(e^{tA}-I)u_{0}}, \\[0,1cm]
\displaystyle{y_{1}(t)=\m \int_{0}^{t}x(s)' A\, u(s) ds}, \\
\displaystyle{y_{2}(t)=\m \int_{0}^{t}x(s)' \til{A} \,u(s) ds},
\end{cases}
\eqn

\subsubsection{Maxwell points}

Consider $\g(t)=E^{u_{0},r,0}(t)$, the geodesic associated with the problem \eqref{eq:qqq3} and with initial covector $(u_{0},r,0), r>0$. Let us first show that there exists another geodesic reaching the point $\g(T^{*})$ in time $T^{*}=2\pi/(a_{1}r)$.
Using Arnol'd's terminology, points reached in the same time by more than one geodesic are called Maxwell points. 
 At the end of this section we prove that $\g$ cannot be optimal after $T^{*}$.
 
Set $u_{0}=(u_{1},u_{2},u_{3},\ldots,u_{m})$
and consider the following variation of the horizontal covector
$$u_{0}^{\omega}=(\cos \omega\, u_{1}+\sin \omega\, u_{2},-\sin \omega\, u_{1}+\cos\omega\, u_{2},u_{3},\ldots,u_{m}), \qquad \omega \in [0,2\pi].$$
Denote $\g^{\omega}(t)=(x^{\omega},y_{1}^{\omega},y_{2}^{\omega}):=E^{u_{0}^{\omega},r,0}(t)$ the geodesic associated with this variation. \\

{\bf Claim}: There exists $\omega\neq 0$ such that $\g(T^{*})=\g^{\omega}(T^{*})$.

\begin{proof}[Proof of the Claim]

Denote by $M_{A}(t):=A^{-1}(e^{tA}-I)$
and notice that
\begin{align*}
M_{A}(t)=A^{-1}(e^{tA}-I)&=
\begin{pmatrix}
\frac{\sin a_{1}rt}{a_{1}r}&\frac{1-\cos a_{1}rt}{a_{1}r}&&&\\
\frac{-1+\cos a_{1}rt}{a_{1}r}&\frac{\sin a_{1}rt}{a_{1}r}&&&\\
&&\ddots&&\\
&&&\frac{\sin a_{\ell}rt}{a_{\ell}r}&\frac{1-\cos a_{\ell}rt}{a_{\ell}r}\\
&&&\frac{-1+\cos a_{\ell}rt}{a_{\ell}r}&\frac{\sin a_{\ell}rt}{a_{\ell}r}
\end{pmatrix}.
\end{align*}
In other words we can write
\begin{align*}
M_{A}(t)&=
\begin{pmatrix}
D_{1}(t)&&\\
&\ddots&\\
&&D_{\ell}(t)\\
\end{pmatrix},
\end{align*}
where
\bqn \label{eq:sincos}
D_{i}(t)=
\begin{pmatrix}
\frac{\sin a_{i}rt}{a_{i}r}&\frac{1-\cos a_{i}rt}{a_{i}r}\\
\frac{-1+\cos a_{i}rt}{a_{i}r}&\frac{\sin a_{i}rt}{a_{i}r}\\
\end{pmatrix}=
2\frac{\sin(a_{i}rt/2)}{a_{i}r}
\begin{pmatrix}
\cos (a_{i}rt/2)&\sin (a_{i}rt/2)\\
-\sin (a_{i}rt/2)&\cos (a_{i}rt/2)\\
\end{pmatrix}.
\eqn
We prove our claim by steps. 

\noindent
(i). From \eqref{eq:sincos} it is easy to see that
$$x^{\omega}(T^{*})-x(T^{*})=A^{-1}(e^{T^{*}A}-I)(u^{\omega}_{0}-u_{0})=M_{A}(T^{*})(u^{\omega}_{0}-u_{0})=0, \qquad \all \omega\in [0,2\pi],$$
since $e^{T^{*}A}-I$ (and so $M_{A}(T^{*})$) has its first $2\times2$ block equal to zero.\\

\noindent
(ii). Now we show that $y^{\omega}_{1}(T^{*})=y_{1}(T^{*})$ for all $\omega\in[0,2\pi]$. Indeed from \eqref{eq:qqq3}, we get
\begin{align*}
y_{1}(t)&=-\m u_{0}'  \int_{0}^{t} (e^{-sA}-I)e^{sA} ds\, u_{0}\\
&=\m u_{0}'  \int_{0}^{t} (e^{sA}-I)ds\,  u_{0}\\
&=\m\la M_{A}(t)u_{0},u_{0}\ra+\m t\|u_{0}\|^{2},
\end{align*}
and
\begin{align*}
y^{\omega}_{1}(T^{*})-y_{1}(T^{*})&=\m\left(\la M_{A}(T^{*})u^{\omega}_{0},u^{\omega}_{0}\ra-\la M_{A}(T^{*})u_{0},u_{0}\ra\right)+\frac{T^{*}}{2}(\|u^{\omega}_{0}\|^{2}-\|u_{0}\|^{2}).
\end{align*}
First notice that
$$\|u^{\omega}_{0}\|^{2}=\|u_{0}\|^{2}, \qquad \all \omega \in[0,2\pi].$$
Moreover, setting $u^{\omega}_{0}=u_{0}+v^{\omega}$ we get (we omit $T^{*}$ in the argument of $M_{A}$)
\begin{align*}
\la M_{A}u^{\omega}_{0},u^{\omega}_{0}\ra-\la M_{A}u_{0},u_{0}\ra
&=\la M_{A}v^{\omega},v^{\omega}\ra+\la (M_{A}+M_{A}')v^{\omega},u_{0}\ra=0
\end{align*}
since the first $2\times2$ block of $M_{A}$ that is zero at $T^{*}$ and $v^{\omega}$ has nonzero component only in the first two entries.

\brem Note that (i) and (ii) are just the manifestation of the fact that, forgetting about the second vertical component $y_{2}$, we are facing the corank 1 case, for which $T^{*}$ is a cut time and there is a rotational symmetry that implies that it is also a conjugate time.
\erem

\noindent
(iii). Now one can proceed in a similar way and compute
\begin{align*}
y_{2}(t)&=\m\int_{0}^{t}x(s)' \til{A}\, u(s) ds\\
&=\m u_{0}' \int_{0}^{t}(A^{-1}(e^{sA}-I))'\til{A}\, e^{sA}ds\,  u_{0}\\
&=-\m u_{0}' \int_{0}^{t}(e^{-sA}-I)A^{-1}\til{A}\,e^{sA}ds\,  u_{0}\\
&=\la C(t) u_{0},u_{0}\ra,
\end{align*}
where we set
\bqn \label{eq:ccc}
C(t)=\m\int_{0}^{t}(e^{-sA}-I)A^{-1}\til{A}\,e^{sA}ds.
\eqn
Since all matrices appearing in \eqref{eq:ccc} but $\til{A}$ are $2\times2$ block diagonal, the first $2\times2$ diagonal block of 
$$K(s):=(e^{-sA}-I)A^{-1}\til{A}\,e^{sA},$$
 is the product of the respective blocks. A direct computation shows that it is
 $$\frac{\z_{0}}{a_{1}}
 \begin{pmatrix}
1-\cos(a_{1}rs)&-\sin(a_{1}rs)\\
\sin(a_{1}rs)&1-\cos(a_{1}rs)\\
\end{pmatrix},
 $$
 where $\begin{pmatrix}
0&\z_{0}\\
-\z_{0}&0\\
\end{pmatrix}$ denotes the first $2\times2$ block of $\til{A}$. Integrating from $0$ to $T^{*}$ one obtains for the first block of $C(T^{*})$
\bqn \label{eq:Cdiag}
\frac{\pi \z_{0}}{a^{2}r^{2}}
 \begin{pmatrix}
1&0\\
0&1
\end{pmatrix}.
\eqn

 

As before, we set $u^{\omega}_{0}=u_{0}+v^{\omega}$ and we get (omiting $T^{*}$ in the argument of $C$)
 \begin{align} \label{eq:hotel}
y^{\omega}_{2}(T^{*})-y_{2}(T^{*})&=\la C v^{\omega},v^{\omega}\ra+\la (C+C')v^{\omega},u_{0}\ra.
\end{align}
Using \eqref{eq:Cdiag} and 
\begin{align*}
\|v^{\omega}\|^{2}&=((\cos \omega- 1) u_{1} + 
     \sin\omega \, u_{2})^2 + (-\sin\omega\, u_{1} + (\cos\omega - 1) u_{2})^2\\
     &=4 (u_{1}^2 + u_{2}^2) \sin^{2}(\omega/2),
\end{align*}
one gets that \eqref{eq:hotel} is linear with respect to the variables
$$\cos \omega -1= -2\sin^{2}(\omega/2), \qquad \sin \omega = 2 \cos (\omega/2) \sin (\omega/2).$$
In other words, if we prescribe that the expression \eqref{eq:hotel} is zero, we get 
\bqn \label{eq:omega}
C_{0}\sin(\omega/2)(C_{1}\cos(\omega/2)+C_{2}\sin (\omega/2))=0,
\eqn
for some suitable constants $C_{0},C_{1},C_{2}$ that do not depend on $\omega$. The Claim is proved since equation \eqref{eq:omega} has always a nontrivial solution $\til{\omega}\in [0,2\pi]$.
\end{proof}
Let us now show that $\g(t)$ cannot be optimal after $T^{*}$. From the previous computation we have $\dot{\g}(T^{*})\neq \dot{\g}^{\til{\omega}}(T^{*})$. By contradiction if $\g$ is optimal after time $T^{*}$ then the concatenation of $\g^{\til{\omega}}|_{[0,T^{*}]}$ and $\g|_{[T^{*},T^{*}+\eps]}$ (for some $\eps>0$) is optimal as well, which is impossible since all optimal trajectories are projections of the Hamiltonian system associated with \eqref{eq:Ham1} and they are smooth.

\subsubsection{Optimality of geodesics}
In this section we prove that $\g(t)=E^{u_{0},r,0}(t)$, $r>0$, is optimal up to its first Maxwell time $T^{*}=2\pi/(a_{1} r)$.\\
To this extent, 
consider the following auxiliary optimal control problem:\\

\noindent
{\bf P.} Let $T<T^{*}$ and set $(\bar{x},\bar{y}_{1},\bar{y}_{2})=\g(T)$. Find a length-parametrized trajectory of the system 
\bqn \label{eq:controlsyst1}
\begin{cases}
\dot{x}_{i}=u_{i}, \qquad\qquad i=1,\ldots,m,\\
\dot{y}_{1}=\m x^{\prime}L_{1}u,\\
\dot{y}_{2}=\m x^{\prime}L_{2}u,
\end{cases}
\eqn
 starting from the origin, and reaching the hyperplane $\{x=\bar{x}\}$ in time $T$, maximizing the $y_{1}$ coordinate.
\brem \label{r:r} Notice that $\bar{y}_{1}>0$ since $r>0$ implies that the trajectory is not contained in the hyperplane $\{y_{1}=0\}$.
\erem
\bl The following assertions hold: (i) There exists a solution $\g^{*}$ of the problem {\bf P}. (ii) $\g^{*}$ is a length minimizer. (iii) $\g^{*}(t)=E^{\til{u}_{0},\til{r},0}(t)$ for some $(\til{u}_{0},\til{r})$.
\el
\begin{proof} Let us prove (i). In problem {\bf P}, since we deal with length-parametrized trajectories, we can assume that the set of controls in \eqref{eq:controlsyst1} is $U=\{\|u\|\leq 1\}$. The existence of a solution of {\bf P} can be obtained with standard arguments using the compactness and convexity of the set of admissible velocities (see \cite{agrachevbook,cesari}).

To prove (ii) assume  by contradiction, 
that there exists a trajectory of \eqref{eq:controlsyst1} reaching the point $(\bar{x},y_{1}^{*},y_{2}^{*})=\g^{*}(T)$ 
in time $T_{0}<T$. 
By small time controllability there exists a trajectory  of system \eqref{eq:controlsyst1} reaching in time $T$ the point $(\bar{x},\wh{y}_{1},\wh{y}_{2})$, with $\wh{y}_{1}>y^{*}_{1}$ contradicting the fact that $\g^{*}$ maximize the $y_{1}$ coordinate. The fact that $\g^{*}$ is also a length minimizer follows from
Remark \ref{r:time}.

To prove (iii) observe that $\g^{*}$ satisfies the Pontryagin Maximum Principle (see again \cite{agrachevbook}) for the problem of minimizing $-y_{1}=-\int_{0}^{T} \dot{y}_{1} dt = -\int_{0}^{T} x'L_{1}u \, dt$, i.e. with the Hamiltonian
\begin{align}\label{eq:delirio}
\mc{H}_{u}&= \sum_{i=1}^{m} \la \lam, u_{i}X_{i}\ra + \nu x'L_{1}u\\
&=p u+r_{1}x'L_{1}u+r_{2}x'L_{2}u+ \nu x'L_{1}u. \nn
\end{align}
where $\lam=(p,r_{1},r_{2})$ are the dual variables to $(x,y_{1},y_{2})$ in $T^{*}M$. In formula \eqref{eq:delirio} $\nu$ is a nonnegative constant. The Hamiltonian equations give 
$$
\begin{cases}
\dot{r}_{1}=-\frac{\partial \mc{H}_{u}}{\partial y_{1}}=0, \\[0.2cm]
\dot{r}_{2}=-\frac{\partial \mc{H}_{u}}{\partial y_{2}}=0,\\[0.2cm]
\dot{p}'=-\frac{\partial \mc{H}_{u}}{\partial x}=-(r_{1}L_{1}+r_{2}L_{2}+ \nu L_{1})u.
\end{cases}
$$
Since the final point is constrained on the set $\{x=\bar{x}\}$, the transversality conditions give $r_{1}=0,r_{2}=0$.
Hence we have
\begin{gather} 
\mc{H}_{u}=(p +\nu x'L_{1})u,  \nn \\
\dot{p}'=-\nu L_{1}u. \label{eq:dd1}
\end{gather}
Notice that actually $\nu>0$, otherwise the trajectory is a straight line contained in the plane $\{y_{1}=y_{2}=0\}$, see Remark \ref{r:r}.
The maximality condition and the condition that the final time is fixed in such a way that trajectories are parameterized by length give
$$\mc{H}_{u(t)}(x(t),y_{1}(t),y_{2}(t),p(t),r_{1},r_{2})=\max_{v}\mc{H}_{v}(x(t),y_{1}(t),y_{2}(t),p(t),r_{1},r_{2})=1,$$
\bqn \label{eq:dd2}
u(t)=\frac{p'(t)-\nu L_{1}x(t)}{\|p'(t)-\nu L_{1}x(t)\|}=p'(t)-\nu L_{1}x(t).
\eqn
Notice that a geodesic for the problem \eqref{eq:controlsyst1} associated with the covector $(u_{0},r_{1},r_{2})$ corresponds to a control 
\bqn \label{eq:dd3}u(t)=p'-r_{1}L_{1}x-r_{2} L_{2}x,\eqn
where
\bqn \label{eq:dd4}\dot{p}'=(-r_{1}L_{1}-r_{2}L_{2})(p'-r_{1}L_{1}x-r_{2}L_{2}x).\eqn
Comparing equations \eqref{eq:dd1} - \eqref{eq:dd2} with \eqref{eq:dd3} - \eqref{eq:dd4} it follows that $\g^{*}$ is a geodesic for the problem \eqref{eq:controlsyst1} corresponding to an initial covector $(\til{u}_{0},\nu,0)$, for some $\til{u}_{0}$, with $\|\til{u}_{0}\|=1$. Then (iii) is proved for $\til{r}=\nu$.
\end{proof}

\noindent
We have the following \\

{\bf Claim}. $\g^{*}= \g$, i.e. $\til{u}_{0}=u_{0}$ and $\til{r}=r$.

\begin{proof}[Proof of the Claim]
It is enough to prove that the parameters $u_{0},r$ such that a geodesic $\g(t)=(x(t),y_{1}(t),y_{2}(t))=E^{u_{0},r,0}$ satisfies $x(T)=\bar{x}$ with $T<T^{*}$ are unique. 

From the computations in Sections \ref{s:cut} we know that
 $$x(t)=M_{A}(t)u_{0}, \qquad \text{where} \quad A=r L_{1}, \quad\text{and}\quad M_{A}(t)=A^{-1}(e^{tA}-I).$$
In particular, using the non singularity of $A$, the equality at $t=T$ gives
\bqn \label{eq:u0} 
u_{0}=M_{A}^{-1}(T)\bar{x}.
\eqn
Computing the norm of vectors in equality \eqref{eq:u0}, it follows
\bqn \label{eq:xyz}
1=\|u_{0}\|^{2}=\sum_{i=1}^{\ell} \frac{\rho_{i}(\bar{x})^{2}}{T^2}\frac{a_{i}rT/2}{\sin(a_{i}rT/2)}, \qquad\text{where}\quad \rho_{i}(x)=(x_{2i-1}^{2}+x_{2i}^{2})^{1/2}.
\eqn
Notice that the right hand side of \eqref{eq:xyz} is the sum of monotonic functions with respect to the variable $rT$, on the segment $[0,2\pi/a_{1}]$  ($T<T^{*}$ implies $rT\leq 2\pi/a_{1}$). 

Moreover since the curve is length-parametrized we have $\|x(T)\| \leq T$. As a consequence there exists a unique solution $rT$ of equation \eqref{eq:xyz} in the segment $[0, 2\pi/a_{1}]$. In particular $r$ is uniquely determined and $u_{0}$ is uniquely recovered from equation \eqref{eq:u0}.
 \end{proof}
Since $\g=\g^{*}$ and $\g^{*}$ is length-minimizer for every $T<T^{*}$, it follows that $\tcut=T^{*}$.

\subsection{First conjugate time} \label{s:conj}
In this section we prove that in the corank 2 case, the cut time is not equal to the first conjugate time, in general.
This is deeply different  from the corank 1 case, where the cut locus always coincides with the first conjugate locus.

It is enough to show that the cut time is not conjugate in the $(4,6)$ case.
Define the Jacobian of the exponential map
\bqn\label{eq:tconj0}
J_{\Exp}(t,u_{0},r_{1},r_{2}):=\det \left(\frac{\partial \Exp}{\partial t},\frac{\partial \Exp}{\partial u_{0}},\frac{\partial \Exp}{\partial r_{1}},\frac{\partial \Exp}{\partial r_{2}}\right).
\eqn
\brem
Recall that the first conjugate time $\tconj$ for the geodesic corresponding to the covector $(u_{0},r_{1},r_{2})$ is the first time $t>0$ for which we have 
\bqn\label{eq:tconj}
J_{\Exp}(t,u_{0},r_{1},r_{2})=0.
\eqn
\erem
\noindent
We have to prove that equation \eqref{eq:tconj} is not satisfied when $t=\tcut$.

To compute $J_{\Exp}$ we use the following trick.
Let $\lv=p dx+r_{1}dy_{1}+r_{2}dy_{2}$ be the Liouville form. 
\bl \label{l:lv} We have 
$
\lv\left(\dfrac{\partial \Exp}{\partial t}\right)=1, \, \lv\left(\dfrac{\partial \Exp}{\partial u_{0}}\right)=\lv\left(\dfrac{\partial \Exp}{\partial r_{i}}\right)=0.
$
\el
\begin{proof} The first equality follows from the fact that the Hamiltonian is homogeneous of degree 2. Indeed set $\lam=(p,r_{1},r_{2})$ and $q=(x,y_{1},y_{2})$, we have $\lv=\lam dq$ and
\bqn \nn
\lv\left(\frac{\partial \Exp}{\partial t}\right)=\la\lam, \frac{\partial q}{\partial t}\ra= \la\lam ,\frac{\partial H}{\partial \lam}\ra=2H=1,
\eqn
since length-parametrized trajectory belong to the set $\{H=1/2\}$. The second and the third identities follow from the fact that the Liouville form is preserved by the Hamiltonian flow, hence the values of $\lv(\frac{\partial \Exp}{\partial u_{0}})$ and $\lv(\frac{\partial \Exp}{\partial r_{i}})$ are constant with respect to $t$. In particular at $t=0$ 
they are annihilated by the Liouville form.
\end{proof}
If we compute the exponential map in a neighborhood of a geodesic with $(r_{1},r_{2})=(r,0)$, with $r\neq 0$, using the identity
$r_{1}dy_{1}=\lv-pdx-r_{2}dx_{2}$ and Lemma \ref{l:lv} we get
\begin{align}
J_{\Exp}&=dx\wedge dy_{1}\wedge dy_{2} \left(\frac{\partial \Exp}{\partial t},\frac{\partial \Exp}{\partial u_{0}},\frac{\partial \Exp}{\partial r_{1}},\frac{\partial \Exp}{\partial r_{2}} \right)\nn \\
&=\frac{1}{r_{1}}dx \wedge dy_{2} \left(\frac{\partial \til{\Exp}}{\partial u_{0}},\frac{\partial \til{\Exp}}{\partial r_{1}},\frac{\partial \til{\Exp}}{\partial r_{2}}\right),\label{eq:Je}
\end{align}
where $\til{\Exp}(t,u_{0},r_{1},r_{2})=(x(t,u_{0},r_{1},r_{2}),y_{2}(t,u_{0},r_{1},r_{2}))$ denote the exponential map  where $y_{1}$ is removed.
More precisely, \eqref{eq:Je} is the function of $(t,u_{0},r_{1},r_{2})$ given by
$$J_{\Exp}=\frac{1}{r_{1}}\det\begin{pmatrix}
\dfrac{\partial x}{\partial u_{0}}v_{1}&\dfrac{\partial x}{\partial u_{0}}v_{2}&\dfrac{\partial x}{\partial u_{0}}v_{3}& \dfrac{\partial x}{\partial r_{1}}& \dfrac{\partial x}{\partial r_{2}} \\[0,4cm]
\dfrac{\partial y_{2}}{\partial u_{0}}v_{1}&\dfrac{\partial y_{2}}{\partial u_{0}}v_{2}&\dfrac{\partial y_{2}}{\partial u_{0}}v_{3} & \dfrac{\partial y_{2}}{\partial r_{1}}& \dfrac{\partial y_{2}}{\partial r_{2}} \\
\end{pmatrix},$$
where $v_{1},v_{2},v_{3}$ are 3 independent tangent vectors to the 3-sphere $\{u_{0}\in \R^{4}, \|u_{0}\|=1\}$. We select
\bqn
v_{1}=\begin{pmatrix} -u_{2}\\u_{1}\\0\\0 \end{pmatrix},\quad
v_{2}=\begin{pmatrix} 0\\0\\-u_{4}\\u_{3} \end{pmatrix},\quad 
v_{3}=\begin{pmatrix} -u_{4}\\0\\0\\u_{1} \end{pmatrix}.
\eqn

From the computation of Section \ref{s:cut} one easily gets
\begin{align*}
\frac{\partial x}{\partial u_{0}}=M_{A}(t)=A^{-1}(e^{tA}-I),
\end{align*}
\begin{align*}
\frac{\partial x}{\partial r_{1}}&=-A^{-1}A(A^{-1}(e^{tA}-I)+t\,e^{tA})u_{0}\\
&=-(M_{A}(t)+t\,e^{tA})u_{0},
\end{align*}
\begin{align*}
\frac{\partial x}{\partial r_{2}}&=-A^{-1}\til{A}(A^{-1}(e^{tA}-I)+t\,e^{tA})u_{0}\\
&=-A^{-1} \til{A}(M_{A}(t)+t\,e^{tA})u_{0}.
\end{align*}
Moreover (see again Section \ref{s:cut})
\begin{align}\label{eq:y2}
y_{2}(t)&
=\la C(t) u_{0},u_{0}\ra,
\end{align}
where
\bqn \label{eq:defC}
C(t)=-\m \int_{0}^{t}(e^{-sA}-I)A^{-1}\til{A}\,e^{sA}ds.
\eqn
The function $y_{2}(t)$ from \eqref{eq:y2}, being a quadratic form with respect to $u_{0}$, gives
$$\dfrac{\partial y_{2}}{\partial u_{0}}v_{i}=\la (C(t)+C'(t))u_{0},v_{i}\ra.$$
Now we compute these derivatives at $t=\tcut=\frac{2\pi}{ar}$ where $a>b$ are the moduli of the eigenvalues of $A$.

It is easily seen that
$$
B:=M_{A}(\tcut)=
\begin{pmatrix}
0&0&0&0\\
0&0&0&0\\[0.3cm]
0&0&\frac{\sin(2\pi b/a)}{br}&\frac{2\sin^{2}(\pi b/a)}{br}\\[0.3cm]
0&0&-\frac{2\sin^{2}(\pi b/a)}{br}&\frac{\sin(2\pi b/a)}{br}
\end{pmatrix},$$
from which it follows that 
$$\begin{pmatrix}
Bv_{1} & Bv_{2} & B v_{3}
\end{pmatrix}=
\begin{pmatrix}
0&0&0\\
0&0&0\\[0.4cm]
0&u_{3}\frac{2\sin^{2}(\pi b/a)}{br}-u4\frac{\sin(2\pi b/a)}{br}&u_{1}\frac{2\sin^{2}(\pi b/a)}{br}\\[0.4cm]
0&u_{4}\frac{2\sin^{2}(\pi b/a)}{br}+u_{3}\frac{\sin(2\pi b/a)}{br}&u_{1}\frac{\sin(2\pi b/a)}{br}\\
\end{pmatrix}
=\begin{pmatrix}
0&0\\

0&M
\end{pmatrix},
$$
where the last identity defines the matrix $M$.

The Jacobian determinant of the exponential map computed at $t=\tcut$ is then expressed as follows
\begin{align}J_{\Exp}(\tcut)&=\det\begin{pmatrix}
\dfrac{\partial x}{\partial u_{0}}v_{1}&\dfrac{\partial x}{\partial u_{0}}v_{2}&\dfrac{\partial x}{\partial u_{0}}v_{3}& \dfrac{\partial x}{\partial r_{1}}& \dfrac{\partial x}{\partial r_{2}} \\[0.4cm]
\dfrac{\partial y_{2}}{\partial u_{0}}v_{1}&\dfrac{\partial y_{2}}{\partial u_{0}}v_{2}&\dfrac{\partial y_{2}}{\partial u_{0}}v_{3} & \dfrac{\partial y_{2}}{\partial r_{1}}& \dfrac{\partial y_{2}}{\partial r_{2}} \\
\end{pmatrix} \nn \\[0.4cm]&=\det
\begin{pmatrix}
\vdots&\vdots&\vdots&\vdots&\vdots\\
Bv_{1}&Bv_{2}&Bv_{3}& \dfrac{\partial x}{\partial r_{1}}& \dfrac{\partial x}{\partial r_{2}} \\
\vdots&\vdots&\vdots&\vdots&\vdots\\[0.4cm]
\dfrac{\partial y_{2}}{\partial u_{0}}v_{1}&\dfrac{\partial y_{2}}{\partial u_{0}}v_{2}&\dfrac{\partial y_{2}}{\partial u_{0}}v_{3} & \dfrac{\partial y_{2}}{\partial r_{1}}& \dfrac{\partial y_{2}}{\partial r_{2}} \\
\end{pmatrix}\nn\\[0.4cm]
&=\det
\begin{pmatrix}
0&0&0&\frac{\partial x_{1}}{\partial r_{1}}& \frac{\partial x_{1}}{\partial r_{2}}\\[0,2cm]
0&0&0&\frac{\partial x_{2}}{\partial r_{1}}& \frac{\partial x_{2}}{\partial r_{2}}\\[0,2cm]
0&M_{11}&M_{12}&\frac{\partial x_{3}}{\partial r_{1}}& \frac{\partial x_{3}}{\partial r_{2}}\\[0,2cm] 
0&M_{21}&M_{22}&\frac{\partial x_{4}}{\partial r_{1}}& \frac{\partial x_{4}}{\partial r_{2}}\\[0,2cm] 
\la Cu_{0},v_{1}\ra&\la Cu_{0},v_{2}\ra&\la Cu_{0},v_{3}\ra& \frac{\partial y_{2}}{\partial r_{1}}& \frac{\partial y_{2}}{\partial r_{2}} \\
\end{pmatrix},\label{eq:conj}
\end{align}
where we use the notation
 $$M=\begin{pmatrix}
M_{11}&M_{12}\\
M_{21}&M_{22}
\end{pmatrix}.$$
From \eqref{eq:conj} it follows
\bqn \label{eq:det0}
J_{\Exp}(\tcut)=\la Cu_{0},v_{1}\ra \cdot \det M \cdot \det N,
\eqn
where $N$ is the matrix
$$N
=\begin{pmatrix}
\frac{\partial x_{1}}{\partial r_{1}}&\frac{\partial x_{1}}{\partial r_{2}}\\[0.2cm]
\frac{\partial x_{2}}{\partial r_{1}}&\frac{\partial x_{2}}{\partial r_{2}}\\
\end{pmatrix}.
$$

It is easy to see from the explicit expression of the geodesics that, in the general case when $a\neq b$, the three factors in \eqref{eq:det0} do not vanish identically in $u_{0}$, since the matrix $\til{A}$ is arbitrary.  The proof of Theorem \ref{t-synthesis} is then completed. 
\brem Notice that $M$ is the zero matrix when $a = b$. Hence, in the $(4,6)$ case, $\tcut = \tconj$ for those $\theta$ such that $L_{\theta}$ has double eigenvalue. Moreover in this case
the rank of the Jacobian matrix drops by 2, since the first three columns are proportional.
\erem

\section{The nilpotent $(4,6)$ case}\label{s-46}
In this section we restrict to the (4,6) case. 
By the previous discussion the geodesics of the sub-Riemannian metric can be written as follows
\bqn \label{eq:qqqtris}
\begin{cases}\displaystyle{
x(t)=A^{-1}(e^{tA}-I)u_{0}} \\[0,1cm]
\displaystyle{y_{1}(t)=\m\int_{0}^{t}x(s)' A\, u(s) ds}, \qquad A=rL_{\theta},\ \til{A}=r\til{L}_{\theta},\\[0,3cm]
\displaystyle{y_{2}(t)=\m \int_{0}^{t}x(s)' \til{A} \,u(s) ds}
\end{cases}
\eqn
and we can assume the matrix $L_{\theta}$ to be diagonal 
\bqn\label{eq:}
L_{\theta}=
\begin{pmatrix}
0&a&0&0\\
-a&0&0&0\\
0&0&0&b\\
0&0&-b&0
\end{pmatrix}, \qquad a\geq b,
\eqn
while $\til{L}_{\theta}$ is an arbitrary skew-symmetric matrix 
$$
\til{L}_{\theta}=
\begin{pmatrix}
0&\z_{0}&\z_{1}&\z_{2}\\
-\z_{0}&0&\z_{3}&\z_{4}\\
-\z_{1}&-\z_{3}&0&\z_{5}\\
-\z_{2}&-\z_{4}&-\z_{5}&0
\end{pmatrix}.$$

\subsection{Proof of Theorem \ref{t-quaternions}}

Recall that the cut time $\tcut$ coincide with $\tconj$ if and only if $\tcut$ is a time that satisfies the equation
\bqn \label{eq:det}
J_{\Exp}(t,u,\theta)\big|_{t=\tcut}=0. 
\eqn

$(P2) \Rightarrow (P1).$ We consider separately the two cases:
\bi
\iii[(a)] $L_{1},L_{2}$ both belong to the same subspace, either $\Q$ or $\Qh$. Then it is not restrictive to assume that $L_{1},L_{2}\in \Q$. In this case all linear combination of $L_{1},L_{2}$ belong to $\Q$, i.e. $L_{\theta} =\cos (\theta) L_{1}+\sin (\theta) L_{2} \in \Q$ for every $\theta$. In particular 
$L_{\theta}$ has a double eigenvalue for every $\theta \in [0,2\pi]$. From the computation of Section \ref{s:cutconj} it is easily seen that $a=b$ implies $M=0$, hence from \eqref{eq:det0} it follows that $\tcut=\tconj$.
\iii[(b)] $L_{1}\in \Q$ and $L_{2}\in \Qh$ ($L_{1}$ and $L_{2}$ plays the same role).
By \eqref{eq:qqhcomm} we have $[L_{1},L_{2}]=0$. Let us prove then that this property implies $(P1)$. 

Indeed every two commuting skew-symmetric matrices can be block diagonalized simultaneously in the same basis. Hence we can assume that, choosing an appropriate coordinate system $y_{1},y_{2}$, that both $L_{1},L_{2}$ are diagonal. As a consequence $L_{\theta}$ and $\til{L}_{\theta}$ are also diagonal. Moreover from \eqref{eq:defC} it is easily seen that, if both $L_{\theta}$ and $\til{L}_{\theta}$ are diagonal, $C$ is $2\times2$ block diagonal, with the first block equal to $c I$, for some constant $c$ (see also \eqref{eq:Cdiag}).

In particular it follows that 
$\la C u_{0},v_{1}\ra=0$
and again \eqref{eq:det0} implies $\tcut=\tconj$.
\ei
$(P1) \Rightarrow (P2).$ By assumption the identity
\bqn \label{eq:det=0}
J_{\Exp}(t,u_{0},\theta)\big|_{t=\tcut(\theta)}=\la C(\theta)u_{0} ,v_{1}\ra \cdot \det M(u_{0},\theta) \cdot \det N(u_{0},\theta)=0,
\eqn
holds for every $u_{0}$ (the horizontal part of the initial covector) and every $\theta$. 
Since the exponential map is linear with respect to $u_{0}$ in the $x$-variable, and quadratic with respect to $u_{0}$ in the $y_{i}$-variables, it follows that \eqref{eq:det=0} is an analytic expression of $(u_{0},\theta)$ (it is polynomial with respect to $u_{0}$ and trigonometric in $\theta$). In particular one of the three factors in \eqref{eq:det=0}  must vanish  identically. 

Assume that $\det M(u_{0},\theta)\equiv0$. Then from the explicit expression it is computed that
$$\det M(u_{0},\theta)= \frac{4 u_{1} u_{3}}{b(\theta)^{2}r^{2}} \sin^{2}\left(\pi \frac{b(\theta)}{a(\theta)}\right),$$
and since $a(\theta)\geq b(\theta)$ by assumption, $\det M\equiv0$ implies $a(\theta)=b(\theta)$, for all $\theta$.

From this it easily follows that $L_{\theta}$ has double eigenvalue for all $\theta$, i.e. if we write 
$$L_{\theta}=q(\theta)+\wh{q}(\theta),$$
it follows that one of $\|q(\theta)\|$ and $\|\wh{q}(\theta)\|$ is identically zero (it is not a restriction to assume $\|\wh{q}(\theta)\|\equiv0$). Hence $L_{\theta}\in\Q$ for all $\theta$, that implies in particular that $L_{1},L_{2} \in \Q$. 

It is not restrictive now to assume that $a(\bar{\theta})\neq b(\bar{\theta})$ for some $\bar{\theta}\in[0,2\pi]$. We show that the identities
\bi
\iii[(a)] $\det N(u_{0},\bar{\theta})=0$,
\iii[(b)] $\la C(\bar{\theta})u_{0},v_{1}\ra=0$,
\ei 
both imply that there exists a choice of the coordinates such that $L_{1}\in \Q$, $L_{2}\in \Qh$.
We give details only for case $(b)$, the other one is similar. Considering $(b)$ as an equation in the variables $\z_{1},\z_{2},\z_{3},\z_{4}$ (the non diagonal entries of the matrix $\til{L}_{\theta}$) it is easy to see that the identity $(b)$ can be written as an equation
$$F(\z_{i},u_{i})=0,$$
where $F$ is a quadratic form in the $u_{i}$ whose coefficients depend linearly on $\z_{i}$. Since these equation should be satisfied for all $u_{0}=(u_{1},\ldots,u_{4})$, choosing values
$$u_{0}\in \{(1,0,1,0),(1,0,0,1),(0,1,1,0),(0,1,0,1)\},$$ 
one gets the set of 4 linear equations:
\bqn
\begin{cases}
(a \z_{4}+b \z_{1}) \cos \eta +(a \z_{3}-b \z_{2}) \sin \eta=0,\\
(a \z_{3}-b \z_{2}) \cos \eta -(a \z_{4}+b \z_{1}) \sin \eta=0,\\
(a \z_{2}-b \z_{3}) \cos \eta +(a \z_{1}+b \z_{4}) \sin \eta=0,\\
(a \z_{1}+b \z_{4}) \cos \eta -(a \z_{2}-b \z_{3}) \sin \eta=0,
\end{cases}
\eqn
where we set $\eta=\pi b/a$, and for simplicity of the notation we denote $a=a(\theta_{0}),b=b(\theta_{0})$.

It is easy to show, using the fact that $a\neq b$, that this system has the unique solution
$$\z_{1}=\ldots=\z_{4}=0,$$
which means that $L_{1}$ and $L_{2}$ are both diagonal. Due to this fact they can be written, as pure quaternions (see Appendix \ref{app1}), as a linear combination of $i,\wh{i}$
\bqn
L_{1}= \alpha i + \wh{\alpha}\, \wh{i}, \qquad L_{2}= \beta i + \wh{\beta}\, \wh{i}. 
\eqn
Performing the change of variables
$$\begin{pmatrix}
\til{y}_{1}\\
\til{y}_{2}
\end{pmatrix}=
\begin{pmatrix}
\alpha & \wh{\alpha}\\
\beta & \wh{\beta}
\end{pmatrix}^{-1}
\begin{pmatrix}
y_{1}\\
y_{2}
\end{pmatrix},
$$
 we find a system of coordinates such that $L_{1}=i$, and $L_{2}= \wh{i}$, i.e. that satisfies (P2).

\subsection{Proof of Theorem \ref{t-haus}}

In the paper \cite{corank1} it is proved that, if the sub-Riemannian manifold is regular\footnote{A sub-Riemannian manifold $(M,\distr,\metr)$ is said to be \emph{regular} if, defining the distributions $\distr^{1}:=\distr, \, \distr^{i+1}:=\distr^{i
}+[\distr^{i},\distr]$, the dimension of $\distr^{i}_{q}$ does not depend on the point for all $i\geq 1$.} with Hausdorff dimension $\h$, the Radon-Nikodym derivative of the spherical Hausdorff measure $\mc{S}^{\h}$ with respect to the Popp's measure $\mu$, denoted $f_{\mc{S}\mu}$, is given by the volume of the unit ball in the nilpotent approximation, namely
\bqn \label{eq:density}
f_{\mc{S}\mu}(q)= \frac{2^{\h}}{\wh{\mu}_{q}(\wh{B}_{q})}, \qquad q\in M,
\eqn
where $\wh{\mu}_{q}$ is the Popp's measure defined on the nilpotent approximation $G_{q}$ of the structure at the point $q$. Note that $\wh{\mu}_{q}$ is the left-invariant measure on $G_{q}$ that coincides with the Popp's measure of the original sub-Riemannian metric at the point $q$. 
\brem
Notice that in our $(m,m+2)$ case the structure is automatically regular since, by assumption, the distribution has constant rank and with one bracket we get all the tangent space. 
\erem
\brem
In \cite{corank1} it is proved that $f_{\mc{S}\mu}$ is a continuous function, which is bounded and bounded away from zero, in restriction to compact sets.
\erem

\brem \label{r:expcov} For the analysis of the regularity of  \eqref{eq:density} it is convenient to parametrize the nilpotent unit ball via the exponential map, as a function defined on the whole fiber in the cotangent space. 
In other words we do not restrict to the set $\{\|u_{0}\|=1\}$ and define
for every $\lam_{0}=(u_{0},r)\in \R^{6}$
$$\Exp(\lam_{0})= \pi(e^{\vec{H}}(\lam_{0})),$$
where $H$ is the Hamiltonian defined in \eqref{eq:Ham1} and $e^{t\vec{H}}$ denotes the flow in $T^{*}M$ of the Hamiltonian vector field associated with $H$. Using the homogeneity property $H(c\lambda)=c^{2}H(\lambda), \ \all c>0$, we have that
$$e^{\vec{H}}(s\lam)=e^{s\vec{H}}(\lam), \qquad \all s>0.$$
In other words we can recover the geodesic on the manifold with initial covector $\lam_{0}$ as the image 
of the ray $\{t\lam_{0},t\in[0,1]\}\subset T^{*}_{q_{0}}M$ that joins the origin to $\lam_{0}$.
$$\Exp(t\lam_{0})=\pi(e^{\vec{H}}(t\lam_{0}))= \pi(e^{t\vec{H}}(\lam_{0}))
=\g(t).$$

\erem
Due to the previous analysis and thanks to Remark \ref{r:expcov}, we can express the volume of the unit ball of the nilpotent approximation at a point $p$ as follows

\bqn \label{eq:VVV}
V(p)=\int_{0}^{2\pi} \int_{0}^{\A(\theta,p)}  \int_{B} J_{\Exp}(u,\theta,r,p) du dr d\theta,
\eqn
where $J_{\Exp}$ is the Jacobian of the exponential map starting from $p$, expressed in the new variables, $B=\{u_{0}=(u_{1},\ldots,u_{4}), \|u_{0}\|\leq 1\}$ is the 4-dimensional unit ball and
$$\A(\theta,p)= \frac{2 \pi}{\max \sigma(L_{\theta}(p))}.$$
The problem of the regularity of the function \eqref{eq:density} is then reduced to the regularity of the function
$$p \mapsto V(p),$$
where $p$ is a (6-dimensional) parameter. 
Since the family of sub-Riemannian metrics is smooth with respect to $p$, the exponential map smoothly depends on the parameter $p$. As a consequence the integrand in \eqref{eq:VVV}, being the Jacobian of the exponential map, is a smooth function of its variables. 

In addition, the function $p\mapsto \A(\theta,p)$ is Lipschitz, being the inverse of the \virg{maximum moduli of eigenvalues} function, which is Lipschitz (see \cite{kurdyka}). In particular $\A(\theta,p)$ admits bounded first derivative almost everywhere with respect to $(\theta,p)$.

\bdeff Define the following sets 
\bi
\iii[-] $\Bad$ is the set of $p$ such that $\exists\, \theta$ for which $L_{\theta}(p)$ has a double eigenvalue,
\iii[-] $\Badfin $ is the set of $p$ such that $\exists$ \emph{a finite number} of $\theta$ for which $L_{\theta}(p)$ has a double eigenvalue,
\iii[-] $\Badinf$ is the set of $p$ such that $\all \theta$,  $L_{\theta}(p)$ has a double eigenvalue.
\ei
\edeff
Thanks to Lemma \ref{l:lemmone},  for a generic sub-Riemannian metric, the set of points $p\in \Badinf$ is a union of isolated points. 
 Moreover, due to the expression \eqref{eq:eig} of the eigenvalues in terms of the quaternions given in the Appendix \ref{app1}, the fact that $L_{\theta}$ has a double eigenvalue for all $\theta$ is written as $\|q(\theta)\|=0$ for all $\theta$ (or the same for $\wh{q}$). This condition is equivalent to the equation $\|q(\theta)\|^{2}=0$, that is analytic in $\theta$. In particular this equation, if it is not identically satisfied, has a finite number of solution in $[0,2\pi]$.

\brem
Notice that the expression \eqref{eq:eig} for the eigenvalues, provides a crucial obstruction for the generalization of the result to $m>4$.
\erem
From this it follows that, for a generic sub-Riemannian metric, the set of critical points $\Bad$ is the disjoint union $\Bad=\Badfin \cup \Badinf$.
Moreover the set of points where $p\mapsto A(\cdot,p)$ is not smooth is contained in $\Sigma$.

\medskip
Let us write the volume function, depending on the parameter $p$, as follows 
\bqn \label{eq:svol}
V(p)=\int_{\theta=0}^{2\pi} \int_{r=0}^{\A(\theta,p)} f(\theta,r,p) dr d\theta,
\eqn
where we denote by 
\begin{gather*}
f(\theta,r,p)= \int_{B} J_{\Exp}(u,\theta,r,p) du.
\end{gather*}
Recall that $f$ is smooth as a function of all its variables, while $A(\theta,p)$ is Lipschitz with respect to the parameters $(\theta,p)$. In particular it has bounded derivatives.


We want to prove that $V$ is $C^{1}$ at any point $p_{0}$. To this extent,
let us write
$$V(p)=\int_{0}^{2\pi} \int_{0}^{\A(\theta,p_{0})} f(\theta,r,p) dr d\theta+\int_{0}^{2\pi} \int_{\A(\theta,p_{0})}^{\A(\theta,p)} f(\theta,r,p) dr d\theta.$$
The function
$$p \mapsto \int_{0}^{2\pi} \int_{0}^{\A(\theta,p_{0})} f(\theta,r,p) dr d\theta,$$
is always smooth since it is the integral of a smooth function (with respect to $p$) on a fixed domain. Denote now
\bqn \label{eq:W}
W(p):=\int_{\theta=0}^{2\pi} \int_{\A(\theta,p_{0})}^{\A(\theta,p)} f(\theta,r,p) dr d\theta.
\eqn
We are left to prove that $W$ is $C^{1}$ around $p_{0}$. Notice that, by definition, $W(p_{0})=0$. 

Assume that $p_{0}\notin \Sigma$. Then, since both functions $A$ and $f$ in \eqref{eq:W} are smooth,  $W$ is $C^{1}$ at $p_{0}$ and the derivative at a point $p$ (in a neighborhood of $p_{0}$), is computed as follows
\bqn\label{eq:volume0}
\frac{ \partial W}{\partial p_{i}}(p)=\int_{0}^{2\pi} \int_{\A(\theta,p_{0})}^{\A(\theta,p)} \frac{\partial f}{\partial p_{i}}(\theta,r,p) dr d\theta+\int_{0}^{2\pi}  \frac{\partial A}{\partial p_{i}}(\theta,p) f(\theta,A(\theta,p),p) dr d\theta.
\eqn

Assume now that $p_{0}\in \Badfin$. The first term in \eqref{eq:volume0} is continuous. Moreover, since at $p_{0}$ there are only a finite number of $\theta$ such that $L_{\theta}$ has double eigenvalue we have
\bqn \label{eq:leb}
 \frac{\partial A}{\partial p_{i}}(\theta,p) \, \underset{p\to p_{0}}{\longrightarrow} \, \frac{\partial A}{\partial p_{i}}(\theta,p_{0}), \qquad \text{a.e.} \ \theta \in [0,2\pi] .
\eqn
Since $\frac{\partial A}{\partial p_{i}}$ is bounded and $f$ is smooth, by Lebesgue's dominated convergence we have that the second term is also continuous. 

Finally, consider the case when $p_{0}\in \Badinf$.
Since $p_{0}$ is an isolated point,  the partial derivatives are defined and continuous in $N_{p_{0}}\setminus \{p_{0}\}$, where $N_{p_{0}}$ is a neighborhood of $p_{0}$. We claim that
\bqn\label{eq:prelemma}
\frac{ \partial W}{\partial p_{i}}(p) \longrightarrow 0, \qquad \text{when} \quad p\to p_{0}\in \Badinf.
\eqn
Indeed, by definition of $\Badinf$, the cut time $A(\theta,p_{0})$ coincides with the conjugate time, i.e. it satisfies the identity
$$J_{\Exp}(u,\theta,A(\theta,p_{0}),p_{0})=0, \qquad \all \theta \in [0,2\pi].$$
From this it follows that $f(\theta,A(\theta,p),p)\to 0$ for all $\theta$, that easily implies \eqref{eq:prelemma}.
From the fact  that $W$ is  continuous  in  $N_{p_{0}}$, $\con^1$ in $N_{p_{0}}\setminus \{p_{0}\}$ and has partial derivatives tending to zero for $p\to p_0$, it follows that $W$ is $\con^1$ in $N_{p_{0}}$.

\section{Appendix}
\subsection{Quaternions}\label{app1}
\label{ap-quaternions}
The Lie algebra $so(4)$ of $4\times 4$ skew-symmetric matrices is the direct sum 
$$so(4)=\Q\oplus \Qh,$$
where $\Q$ is the space of pure quaternions and $\Qh$ is the set of pure skew quaternions.  

The space $\Q$ (resp. $\Qh$) is generated by the three matrices $i,j,k$ (respectively  $\wh{i},\wh{j},\wh{k}$) 

$$
i=\left(
\begin{array}{cccc}
 0 & -1 & 0 & 0 \\
 1 & 0 & 0 & 0 \\
 0 & 0 & 0 & -1 \\
 0 & 0 & 1 & 0
\end{array}
\right),~~j=\left(
\begin{array}{cccc}
 0 & 0 & -1 & 0 \\
 0 & 0 & 0 & 1 \\
 1 & 0 & 0 & 0 \\
 0 & -1 & 0 & 0
\end{array}
\right), ~~k=\left(
\begin{array}{cccc}
 0 & 0 & 0 & -1 \\
 0 & 0 & -1 & 0 \\
 0 & 1 & 0 & 0 \\
 1 & 0 & 0 & 0
\end{array}
\right),
$$
and 
$$
\wh{i}=\left(
\begin{array}{cccc}
 0 & -1 & 0 & 0 \\
 1 & 0 & 0 & 0 \\
 0 & 0 & 0 & 1 \\
 0 & 0 & -1 & 0
\end{array}
\right),~~\wh{j}=\left(
\begin{array}{cccc}
 0 & 0 & 1 & 0 \\
 0 & 0 & 0 & 1 \\
 -1 & 0 & 0 & 0 \\
 0 & -1 & 0 & 0
\end{array}
\right), ~~\wh{k}=\left(
\begin{array}{cccc}
 0 & 0 & 0 & -1 \\
 0 & 0 & 1 & 0 \\
 0 & -1 & 0 & 0 \\
 1 & 0 & 0 & 0
\end{array}
\right).
$$

If we endow $so(4)$ with the Hilbert-Schmitd scalar product
$$\langle L_1,L_2\rangle=\frac{1}{4}\tr(L'_1L_2),$$
then $i,j,k,\wh{i},\wh{j},\wh{k}$ is an orthonormal basis. 

The eigenvalues $\omega_1,\omega_2$  of $A=q+\wh{q}$ satisfy:
\bqn \label{eq:eig}
-(\omega_{1,2})^2
=(\| q\|\pm\|\wh{q}\|)^2.
\eqn
As a consequence an element $A\in so(4)$ has a double eigenvalue if and only if $A\in\Q\cup\Qh$.

Also pure quaternions and pure skew quaternions commute:
\bqn
 \label{eq:qqhcomm}
[q,\wh{q}]=0, ~~q\in \Q,~~\wh{q}\in\Qh.
\eqn

\subsection{Transersality lemma}
\label{ap-transversal}
Let $S$ be the set of $(m,n)$ smooth sub-Riemannian metrics over $M$, equipped with the Whitney topology. Due to the ${\cal C}^\infty$ structure, we have the existence of smooth bump functions and the  results in this section are  essentially local. Then we can assume that $S$ is the set of $m$-tuples  $F=(X_1,\ldots X_m)$ of smooth independent vector fields on some open subset $M$ of $\R^n$, satisfying
$$
T_qM=\Delta_q+[\Delta,\Delta]_q,~~\mbox{ for every }q\in M.
$$

The vector fields $X_1,\ldots,X_m$ form an orthonormal basis for the sub-Riemannian metric $\metr$ they specify.  

Let $B$ be the  bundle over $M$ whose fiber at $q\in M$ is the variety of $k$-dimensional vector space of  $\metr$-skew symmetric endomorphisms of $\Delta_q$. 

Let us consider the mapping 
$$
\ba{rcl}
\rho:  S\times M&\to& B \\
(F,q)&\to&{\cal L}^F_q
\ea
$$
where ${\cal L}^F_q$ has been defined in Remark \ref{r-LFQ}.
It is clear that $\rho$ is ${\cal C}^\infty$.

Let us fix a point $F_0\in S$, a point $q_0\in M$ and coordinates $(x,y)$ in $M$ such that $q_0=(0,0)$ and the nilpotent approximation of $F_0$ reads in control form
\bqn
\label{e-nil}
\begin{cases}
\dot{x}_{i}=u_{i}, \qquad i=1,\ldots,m,\\
\dot{y}_1=\m x'L_1 u,\\
\ \quad\vdots\\
\dot{y}_{k}=\m x'L_{k} u,
\end{cases}
\eqn
In the coordinates $y$, the space ${\cal L}_q^{F_0}$ is the vector subspace spanned by the matrices $L_1,\ldots L_{k}$.
We have a natural gradation in formal power series of $(x,y)$ induced by setting that the $x_i$ have weight 1 and the $y_i$ have weight 2. This induces a formal gradation on formal vector fields on $M_{q_0}$ in which $\frac{\partial}{\partial x_i}$ have weight $-1$  and $\frac{\partial}{\partial y_i}$ have weight $-2$. The vector fields of the nilpotent approximation \r{e-nil} have weight $-1$.

In control form the sub-Riemannian metric $F_0$ itself reads
\bqn
\left(
\ba{c}
\dot x\\
\dot y_{1}\\
\vdots\\
\dot y_{k}
\ea
\right)
=\left(\ba{c}
u\\
\m x'L_1 u\\
\vdots\\
\m x'L_{k} u
\ea
\right)+H
\eqn
where $H$ is a term of order $>-1$ as a formal $u$-dependent vector field.
 Then we take a smooth bump function $b(x,y)$ which is compactly supported in $M$ and which is 1 in a neighborhood of  $q_0=(0,0)$. We consider the affine space ${\cal A}$ of variations of $F_0$ of the form
 
 \bqn
\left(
\ba{c}
\dot x\\
\dot y_{1}\\
\vdots\\
\dot y_{k}
\ea
\right)
=\left(\ba{c}
u\\
\m x'L_1 u\\
\vdots\\
\m x'L_{k} u
\ea
\right)+
\left(
\ba{c}
0\\
\m x' \delta L_1 u\\
\vdots\\
\m x'\delta L_{k} u
\ea
\right)b(x,y)+H.
\eqn
This defines new sub-Riemannian metrics $F_0+\delta F$. Since $S$ is open in the set of all rank $m$ smooth sub-Riemannian metrics over $M$, then, for a small perturbation $\delta F$, we have $F_0+\delta F\in S$.

To show that $\rho_{q_0}$, defined by $\rho_{q_0}(F):=\rho(q_0,F)$, is a submersion at $q_0$ on the fiber $B_{q_0}$, it is enough to observe that 
$$
\ba{rcl}
{\widehat \rho}_{q_0}:{\cal A}&\to&so(4)^{k}\\
\delta F&\mapsto&(L_1+\delta L_1,\ldots,L_{k}+\delta L_{k})
\ea
$$
is an affine submersion. Then we have proven the following Lemma.
\bl
\label{l-trans}
The map $\rho$ is a submersion.
\el
Now let us restrict to the (4,6) case.
\bdeff
We say that a point $q_0\in M$ is critical for a sub-Riemannian metric $F$ if all elements of the subspace $\mc{L}^{F}_{q_0}$ (from Remark \ref{r-LFQ}) have a double eigenvalue.
\edeff
By formula \ref{eq:eig}, this means that, whatever the coordinates $y$, the matrices $L_1(q_0),L_2(q_0)$ both belong either to $\Q$ or to $\Qh$. 

The dimension $d_1$ of 
the fiber of the bundle $B$ is the dimension of the Grassmannian $G(2,6)$ of 2-subspaces of $\R^6$, i.e. $d_1=8$.

The dimension of the set of pairs $L_1,L_2$ that both belong to $\Q$ (respectively $\Qh$), is the dimension $d_2$ of the Grassmannian $G(2,3)$, i.e $d_2=2$.

Let us define now the partially algebraic ``wrong set'' $W\subset B$ as follows: the fiber $W_{q_0}$ is the set of 2-subspaces of the $\metr$-skew symmetric endomorphisms of $\Delta_{q_0}$, whose elements have a double eigenvalue. The codimension of $W$ in $B$ is $d_1-d_2=6$.

The next Lemma follows from Lemma \ref{l-trans} and a non-compact version of Abraham's parametric transversality Theorems (\cite{abraham}).

\bl[$(4,6)$ case] \label{l:lemmone}
The set of sub-Riemannian metrics that have only isolated critical points is open and dense in $S$.
\el

%
%
{\small
\bibliography{bibliogr-corank2-v1}
\bibliographystyle{abbrv}
}
\end{document}